\documentclass[12pt]{amsart}
\newtheorem{theorem}{Theorem}
\newtheorem{lemma}{Lemma}
\newtheorem{proposition}{Proposition}
\begin{document}

\title{}
\author{}
\centerline{\Large Indefinite locally conformal K\"ahler
manifolds} \vskip 1cm \centerline{\large Sorin
Dragomir\footnote{Universit\`a degli Studi della Basilicata,
Dipartimento di Matematica, Campus Macchia Romana, 85100 Potenza,
Italy, e-mail: {\tt dragomir@unibas.it}} \hspace{1cm} Krishan L.
Duggal\footnote{University of Windsor, Windsor, Ontario N9B3P4,
Canada, e-mail: {\tt yq8@uwindsor.ca}}}
\begin{abstract} We study the basic properties of an indefinite
locally conformal K\"ahler (l.c.K.) manifold. Any indefinite
l.c.K. manifold $M$ with a parallel Lee form $\omega$ is shown to
possess two canonical foliations $\mathcal F$ and ${\mathcal
F}_c$, the first of which is given by the Pfaff equation $\omega =
0$ and the second is spanned by the Lee and the anti-Lee vectors
of $M$. We build an indefinite l.c.K. metric on the noncompact
complex manifold $\Omega_+ = (\Lambda_+ \setminus \Lambda_0
)/G_\lambda$ (similar to the Boothby metric on a complex Hopf
manifold) and prove a CR extension result for CR functions on the
leafs of $\mathcal F$ when $M = \Omega_+$ (where $\Lambda_+
\setminus \Lambda_0 \subset {\mathbb C}^n_s$ is $- |z_1 |^2 -
\cdots - |z_s |^2 + |z_{s+1}|^2 + \cdots + |z_n |^2 > 0$). We
study the geometry of the second fundamental form of the leaves of
$\mathcal F$ and ${\mathcal F}_c$. In the degenerate cases
(corresponding to a lightlike Lee vector) we use the technique of
screen distributions and (lightlike) transversal bundles developed
by A. Bejancu et al., \cite{kn:BeDu}.
\end{abstract}
\maketitle

\section{The first canonical foliation}
Let $M$ be a complex $n$-dimensional indefinite Hermitian manifold
of index $0 < \nu < 2n$, with the complex structure $J$ and the
semi-Riemannian metric $g$. As well known $\nu$ must be even, $\nu
= 2s$. $M$ is an {\em indefinite K\"ahler manifold} if $\nabla J =
0$, where $\nabla$ is the Levi-Civita connection of $(M , g)$, cf.
M. Baros \& A. Romero, \cite{kn:BaRo}. An indefinite Hermitian
manifold $M$ is an {\em indefinite locally conformal K\"ahler}
(l.c.K.) {\em manifold} if for any point $x \in M$ there is an
open neighborhood $U$ of $x$ in $M$ and a $C^\infty$ function $f :
U \to {\mathbb R}$ such that $(U , e^{-f} g)$ is an indefinite
K\"ahler manifold.
\par
Note that any two conformally related indefinite K\"ahler metrics
are actually homothetic. Indeed, let $(U , z^1 , \cdots , z^n )$
be a local system of complex coordinates on $M$ and set
$g_{j\overline{k}} = g(Z_j , \overline{Z}_k )$, where $Z_j$ is
short for $\partial /\partial z^j$ (and overbars denote complex
conjugates). If $\hat{g} = e^f g$ then the Levi-Civita connections
$\nabla$ and $\hat{\nabla}$ (of $(M , g)$ and $(M , \hat{g})$,
respectively) are related by
\[ \hat{\nabla}_{Z_j} \overline{Z}_k = \nabla_{Z_j} \overline{Z}_k
- \frac{1}{2} \{ (Z_j f) \overline{Z}_k + (\overline{Z}_k f) Z_j -
g_{j\overline{k}} \nabla f \} \] where $\nabla f$ is the gradient
of $f$ with respect to $g$, i.e. $g(\nabla f , X) = X(f)$, for any
$X \in T(M)$. Let $T(M) \otimes {\mathbb C}$ be the complexified
tangent bundle. Let $Z^{1,0}$ denote the $(1,0)$-component of $Z
\in T(M) \otimes {\mathbb C}$ with respect to the direct sum
decomposition $T(M) \otimes {\mathbb C} = T^{1,0}(M) \oplus
T^{0,1}(M)$, where $T^{1,0}(M)$ is the holomorphic tangent bundle
over $M$ and $T^{0,1}(M) = \overline{T^{1,0}(M)}$. When $g$ and
$\hat{g}$ are indefinite K\"ahler metrics, both $\nabla$ and
$\hat{\nabla}$ descend to connections in $T^{0,1}(M)$ hence
\[ (\overline{Z}_k f) Z_j - g_{j\overline{k}} (\nabla f)^{1,0} = 0 \]
or $\delta_j^s \overline{Z}_k f - g_{j\overline{k}}
g^{s\overline{r}} \overline{Z}_r f = 0$. Contraction of $s$ and
$j$ leads to $(n-1) \overline{Z}_k f = 0$, i.e. $f$ is a real
valued holomorphic function, hence a constant.
\par
Let $(M , J, g)$ be an indefinite l.c.K. manifold. Let then $\{
U_i \}_{i \in I}$ be an open cover of $M$ and $\{ f_i \}_{i \in
I}$ a family of $C^\infty$ functions $f_i : U_i \to {\mathbb R}$
such that $g_i := e^{-f_i} g$ is an indefinite K\"ahler metric on
$U_i$. Then $g_j = e^{f_i - f_j} g_i$ on $U_i \cap U_j$, i.e.
$g_i$ and $g_j$ are conformally related indefinite K\"ahler
metrics on $U_i \cap U_j$, hence $f_i - f_j = c_{ij}$, for some
$c_{ij} \in {\mathbb R}$. In particular $d f_i = d f_j$, i.e. the
local $1$-forms $d f_i$, $i \in I$, glue up to a globally defined
closed $1$-form $\omega$ on $M$. By analogy with the positive
definite case (cf. e.g. S. Dragomir \& L. Ornea, \cite{kn:DrOr})
we shall refer to $\omega$ as the {\em Lee form} of $M$. An
indefinite l.c.K. metric $g$ is called {\em globally conformal
K\"ahler} (g.c.K.) if the Lee form $\omega$ is exact (cf. e.g.
\cite{kn:Vai2} for the Riemannian case). The tangent vector field
$B$ on $M$ defined by $g(X , B) = \omega (X)$, for any $X \in
T(M)$, is the {\em Lee field}. Let us set $c = g(B,B) \in C^\infty
(M)$ and ${\rm Sing}(B) = \{ x \in M : B_x = 0 \}$. Note that
(opposite to the positive definite case) it may be $c = 0$ and
${\rm Sing}(B) = \emptyset$ (when $B$ is lightlike).
\par
Let $\nabla^i$ be the Levi-Civita connection of $(U_i , g_i )$, $i
\in I$. Then
\[ \nabla^i_X Y = \nabla_X Y - \frac{1}{2} \{ X(f_i ) Y + Y(f_i )
X - g(X,Y) \nabla f_i \} \] for any $X, Y \in T(U_i )$, hence the
local connections $\nabla^i$, $i \in I$, glue up to a globally
defined linear connection $D$ on $M$ given by
\[ D_X Y = \nabla_X Y - \frac{1}{2}\{ \omega (X) Y + \omega (Y) X
- g(X,Y) B \} , \] the {\em Weyl connection} of $M$. Clearly $D J
= 0$.
\par
Let us analyze indefinite l.c.K. manifolds with $\nabla \omega =
0$ (the indefinite counterpart of generalized Hopf manifolds, cf.
I. Vaisman, \cite{kn:Vai2}). Such manifolds carry a natural
foliation $\mathcal F$ defined by the Pfaff equation $\omega = 0$.
Also $c \in {\mathbb R}$, so that $B$ is spacelike (respectively
timelike, or lightlike) when $c > 0$ (respectively $c < 0$, or $c
= 0$). We shall prove
\begin{theorem}
Let $M$ be a complex $n$-dimensional indefinite l.c.K. manifold of
index $2s$, $0 < s < n$, with a parallel Lee form and with ${\rm
Sing}(B) = \emptyset$. Then either {\rm i)} $c \neq 0$ and then
each leaf $L$ of $\mathcal F$ is a totally geodesic
semi-Riemannian hypersurface of $(M , g)$ of index
\begin{equation} \label{e:1} {\rm ind}(L) =
\begin{cases} 2s & c
> 0 \cr 2s-1 & c < 0. \cr
\end{cases} \end{equation}
or {\rm ii)} $c = 0$ and then each leaf of $\mathcal F$ is a
totally geodesic lightlike hypersurface of $(M , g)$. \label{t:1}
\end{theorem}
\noindent {\em Proof}. Assume first that $c \neq 0$. Let us show
that \begin{equation} T(M) = T({\mathcal F}) \oplus {\mathbb R} B.
\label{e:2}
\end{equation}
To this end, let $X \in T(M)$. Then $X - \frac{1}{c} \omega (X) B
\in T({\mathcal F})$. Moreover, if $X \in T({\mathcal F}) \cap
{\mathbb R}B$ then $X = \lambda B$, for some $\lambda \in C^\infty
(M)$, and $0 = \omega (X) = \lambda c$ yields $\lambda = 0$, i.e.
$X = 0$. Therefore (\ref{e:2}) holds. Since ${\mathbb R}B$ is
nondegenerate it follows that $T({\mathcal F}) = ({\mathbb
R}B)^\bot$ is nondegenerate, as well, hence each leaf of $\mathcal
F$ is a semi-Riemannian hypersurface of $(M , g)$ and (\ref{e:1})
holds. Let $L$ be a leaf of $\mathcal F$. Let $\nabla^L$ be the
induced connection and $h^L$ the second fundamental form of $i : L
\hookrightarrow M$. If $X, Y \in T(L)$ then, by $\nabla \omega =
0$ and the Gauss formula $\nabla_X Y = \nabla^L_X Y + h^L (X,Y)$
\[ 0 = X(\omega (Y)) = \omega (\nabla_X Y ) = \omega (h^L (X,Y))
\]
hence $h^L = 0$.
\par
Let us assume now that $c = 0$, so that $B \in T({\mathcal F})$.
Let us set as customary (cf. e.g. \cite{kn:BeDu}, p. 140)
\[ {\rm Rad}\; T({\mathcal F})_x = T({\mathcal F})_x \cap
T({\mathcal F})^\bot_x \, , \;\;\; x \in M. \] Clearly $B \in {\rm
Rad} \; T({\mathcal F})$. Note that ${\rm dim}_{\mathbb R} \;
T({\mathcal F})_x = 2n-1$ hence (cf. e.g. Proposition 2.2 in
\cite{kn:BeDu}, p. 6) ${\rm dim}_{\mathbb R} \; T({\mathcal
F})^\bot_x = 1$, for any $x \in M$. Therefore, if ${\rm Sing}(B) =
\emptyset$ then $T({\mathcal F})^\bot = {\mathbb R}B$ and (by
Proposition 1.1 in \cite{kn:BeDu}, p. 78) each leaf of $\mathcal
F$ is a lightlike hypersurface of $(M , g)$.
\par
If $Y \in T({\mathcal F})^\bot$ then $Y$ is orthogonal on the Lee
field, hence $Y \in T({\mathcal F})$. It follows that $T({\mathcal
F})^\bot \subset T({\mathcal F})$, i.e. ${\rm Rad} \; T({\mathcal
F}) = T({\mathcal F})^\bot$. Let $S(T {\mathcal F})$ be a
distribution on $M$ such that
\begin{equation}
\label{e:3}
T({\mathcal F}) = S(T {\mathcal F}) \oplus_{\rm orth}
T({\mathcal F})^\bot .
\end{equation}
If $V$ is a semi-Euclidean space and $W_a \subset V$, $a \in \{ 1
, 2 \}$, are two subspaces then we write $V = W_1 \oplus_{\rm
orth} W_2$ whenever $V = W_1 \oplus W_2$ and the subspaces $W_a$
are mutually orthogonal. According to the terminology in
\cite{kn:BeDu}, p. 78, the portion of $S(T {\mathcal F})$ over a
leaf $L$ of $\mathcal F$ is a {\em screen distribution} on $L$.
The choice of $S(T {\mathcal F})$ is not unique, yet (by
Proposition 2.1 in \cite{kn:BeDu}, p. 5) $S(T {\mathcal F})$ is
nondegenerate, hence
\begin{equation}  T(M) = S(T{\mathcal F}) \oplus_{\rm orth} S(T{\mathcal
F})^\bot . \label{e:4}
\end{equation}
Note that $S(T{\mathcal F})^\bot$ has rank two and $T({\mathcal
F})^\bot \subset S(T{\mathcal F})^\bot$. The following result is
an adaptation (to the foliation $\mathcal F$ of $M$, rather than a
single lightlike hypersurface) of Theorem 1.1 in \cite{kn:BeDu},
p. 79.
\begin{lemma} Let $\pi : E \to M$ be a subbundle of $S(T{\mathcal F})^\bot \to M$
such that $S(T{\mathcal F})^\bot = T({\mathcal F})^\bot \oplus E$.
Let $V \in \Gamma^\infty (U , E)$ be a locally defined nowhere
zero section, defined on the open subset $U \subseteq M$. Then
{\rm i)} $\omega (V) \neq 0$ everywhere on $U$. Let us consider
$N_V \in \Gamma^\infty (U , S(T{\mathcal F})^\bot )$ given by
\begin{equation}
\label{e:5} N_V = \frac{1}{\omega (V)} \left\{ V - \frac{g(V,V)}{2
\omega (V)} \; B \right\} .
\end{equation}
If $V^\prime \in \Gamma^\infty (U^\prime , E)$ is another nowhere
zero section, defined on the open subset $U^\prime \subseteq M$
such that $U \cap U^\prime \neq \emptyset$, then {\rm ii)} $N_V =
N_{V^\prime}$ on $U \cap U^\prime$. Moreover, let $x \in M$ and $U
\subseteq M$ an open neighborhood of $x$ such that $\left.
E\right|_U = \pi^{-1}(U)$ is trivial. Let $V$ us set
\begin{equation} {\rm tr}(T {\mathcal F})_x = {\mathbb R} N_V (x) \, .
\label{e:6}
\end{equation}
Then {\rm iii)} ${\rm tr}(T {\mathcal F})_x$ is well defined and
gives a lightlike subbundle ${\rm tr}(T{\mathcal F}) \to M$ of
$S(T{\mathcal F})^\bot \to M$ such that
\begin{equation}
S(T{\mathcal F})^\bot = T({\mathcal F})^\bot \oplus {\rm tr}(T
{\mathcal F}). \label{e:7}
\end{equation}
Finally, {\rm iv)} the definition of ${\rm tr}(T {\mathcal F})$
doesn't depend upon the choice of complement $E$ to $T({\mathcal
F})^\bot$ in $S(T {\mathcal F})^\bot$. \label{l:1}
\end{lemma}
\noindent {\em Proof}. The proof of (i) is by contradiction. If
$\omega (V)_{x_0} = 0$ for some $x_0 \in U$ then $V_{x_0} \in
T({\mathcal F})_{x_0}$ and then (by (\ref{e:3}))
\[ V_{x_0} \in S(T{\mathcal F})_{x_0} \cap S(T{\mathcal
F})^\bot_{x_0} = (0), \] a contradiction. To prove (ii) let
$V^\prime$ be a nowhere zero section in $E$ on $U^\prime$ (with $U
\cap U^\prime \neq \emptyset$). Then $V^\prime = \alpha V$, for
some $C^\infty$ function $\alpha : U \cap U^\prime \to {\mathbb R}
\setminus \{ 0 \}$, and an inspection of (\ref{e:5}) leads to $N_V
= N_{V^\prime}$ on $U \cap U^\prime$. This also shows that the
definition of ${\rm tr}(T {\mathcal F})_x$ doesn't depend upon the
particular local trivialization chart of $E$ at $x$. To check the
remaining statement in (iii) note first that
\begin{equation}
g(N_V , N_V ) = 0, \;\;\; \omega (N_V ) = 1. \label{e:8}
\end{equation}
The first relation in (\ref{e:8}) shows that ${\rm tr}(T{\mathcal
F})$ is lightlike, while the second relation yields $T({\mathcal
F})^\bot \cap {\rm tr}(T {\mathcal F}) = (0)$. Yet both bundles
have rank one, hence (\ref{e:7}) holds. Finally, if $F \to M$ is
another complement to $T({\mathcal F})^\bot$ in $S(T{\mathcal
F})^\bot$ then an {\em a priori} new lightlike bundle, similar to
${\rm tr}(T{\mathcal F})$, may be build in terms of a local
section $W \in \Gamma^\infty (U , F)$. Yet (by (\ref{e:7})) $W =
\alpha N_V + \beta B$, for some $\alpha , \beta \in C^\infty (U)$,
hence (by (\ref{e:8})) $N_W = N_V$. $\square$ \vskip 0.1in
According to the terminology in \cite{kn:BeDu}, p. 79, the portion
of ${\rm tr}(T{\mathcal F})$ over a leaf $L$ of $\mathcal F$ is
the {\em lightlike transversal vector bundle} of $L$ with respect
to the screen distribution $\left. S(T{\mathcal F})\right|_L$. By
(\ref{e:3})-(\ref{e:4}) and (\ref{e:7}) we obtain the
decomposition
\begin{equation}
T(M) = S(T{\mathcal F}) \oplus_{\rm orth} \left[ T({\mathcal
F})^\bot \oplus {\rm tr}(T{\mathcal F})\right] = T({\mathcal F})
\oplus {\rm tr}(T{\mathcal F}). \label{e:9}
\end{equation}
Let ${\rm tan} : T(M) \to T({\mathcal F})$ and ${\rm tra} : T(M)
\to {\rm tr}(T {\mathcal F})$ be the projections associated with
(\ref{e:9}). Next, we set
\[ \nabla^{\mathcal F}_X Y = {\rm tan} \left( \nabla_X Y \right) ,
\;\;\; h (X,Y) = {\rm tra}\left( \nabla_X Y \right),
\]
\[ A_V X = - {\rm tan}\left( \nabla_X V \right) ,
\;\;\; \nabla^{\rm tr}_X V = {\rm tra} \left( \nabla_X V \right) ,
\]
for any $X, Y \in T({\mathcal F})$ and any $V \in {\rm
tr}(T{\mathcal F})$. Then $\nabla^{\mathcal F}$ is a connection in
$T({\mathcal F}) \to M$, $h$ is a symmetric ${\rm tr}(T{\mathcal
F})$-valued bilinear form on $T({\mathcal F})$, $A_V$ is an
endomorphism of $T({\mathcal F})$, and $\nabla^{\rm tr}$ is a
connection in ${\rm tr}(T{\mathcal F}) \to M$. Also one has
\[ \nabla_X Y = \nabla^{\mathcal F}_X Y + h(X,Y), \;\;\; \nabla_X
V = - A_V X + \nabla^{\rm tr}_X V , \] the {\em Gauss} and {\em
Weingarten formulae} of $\mathcal F$ in $(M , g)$. Clearly, the
pointwise restrictions of $\nabla^{\mathcal F}$, $\nabla^{\rm
tr}$, $h$ and $A_V$ to a leaf $L$ of $\mathcal F$ are respectively
the induced connections, the second fundamental form and the shape
operator of $L$ in $(M , g)$, cf. \cite{kn:BeDu}, p. 83. A leaf
$L$ is {\em totally geodesic} if each geodesic of
$\nabla^{\mathcal F}$ lying on $L$ is also a geodesic of the
semi-Riemannian manifold $(M , g)$.
\par
Let us prove the last statement in Theorem \ref{t:1}. As $\nabla
\omega = 0$ it follows that $\omega (h(X,Y)) = 0$. Yet locally
(with the notations in the proof of Lemma \ref{l:1}) $h(X,Y) =
C(X,Y) N_V$, for some $C(X,Y) \in C^\infty (U)$, hence (by
(\ref{e:8})) $h = 0$ and then by a result in \cite{kn:BeDu1} (cf.
also Theorem 2.2 in \cite{kn:BeDu}, p. 88) each leaf of $\mathcal
F$ is totally geodesic in $(M , g)$. $\square$ \vskip 0.1in We end
this section by the following remark. By Theorem 2.2 in
\cite{kn:BeDu}, p. 88, if $c = 0$ then $\nabla^{\mathcal F}$ is
the Levi-Civita connection of the tangential metric induced by $g$
on $T({\mathcal F})$ and the distribution $T({\mathcal F})^\bot$
is Killing.

\section{Indefinite Hopf manifolds}
Let ${\mathbb C}^n_s$ denote ${\mathbb C}^n$ together with the
real part of the Hermitian form
\[ b_{s,n} (z,w) = - \sum_{j=1}^s z_j \overline{w}_j +
\sum_{j=s+1}^n z_j \overline{w}_j \, , \;\;\; z,w \in {\mathbb
C}^n . \] Let $\Lambda = \{ z \in {\mathbb C}^n \setminus \{ 0 \}
: - \sum_{j=1}^s |z_j |^2 + \sum_{j=s+1}^n |z_j |^2 = 0 \}$ be the
null cone in ${\mathbb C}^n_s$ and $\Lambda_0 = \Lambda \cup \{ 0
\}$. Given $\lambda \in {\mathbb C} \setminus \{ 0 \}$ \[
F_\lambda (z) = \lambda z, \;\;\; z \in {\mathbb C}^n \setminus
\Lambda_0 \, , \] is a holomorphic transformation of ${\mathbb
C}^n \setminus \Lambda_0$. Let $G_\lambda = \{ F_\lambda^m : m \in
{\mathbb Z} \}$ be the discrete group generated by $F_\lambda$.
Then
\begin{theorem} Let $n > 1$, $0 < s < n$ and
$\lambda \in {\mathbb C} \setminus \{ 0 \}$, $|\lambda | \neq 1$.
Then $G_\lambda$ acts freely on ${\mathbb C}^n \setminus
\Lambda_0$ as a properly discontinuous group of holomorphic
transformations, hence the quotient space ${\mathbb C}H^n_s
(\lambda ) = ({\mathbb C}^n \setminus \Lambda_0 )/G_\lambda$ is a
complex manifold and
\begin{equation}
g_{s,n} = |z|_{s,n}^{-2} \left( - \sum_{j=1}^s d z^j \odot d
\overline{z}^j + \sum_{j=s+1}^n d z^j \odot d \overline{z}^j
\right)
\label{e:10}
\end{equation} {\rm (}where $|z|_{s,n} =
|b_{s,n}(z,z)|^{1/2}${\rm )} is a globally defined semi-Riemannian
metric, making ${\mathbb C}H^n_s (\lambda )$ into an indefinite
locally conformal K\"ahler manifold. Moreover, if $0 < \lambda <
1$ then ${\mathbb C}H^n_s (\lambda ) \approx \Sigma^{2n-1} \times
S^1$ {\rm (}a diffeomorphism{\rm )}, where $\Sigma^{2n-1} = \{ z
\in {\mathbb C}^n : |z|_{s,n} = 1 \}$. In particular ${\mathbb
C}H^n_s (\lambda )$ is noncompact. If $\Lambda_+ = \{ z \in
{\mathbb C}^n : b_{s,n}(z,z) \geq 0 \}$ and $\Lambda_- = \{ z \in
{\mathbb C}^n : b_{s,n}(z,z) \leq 0 \}$ {\rm (}so that $\partial
\Lambda_\pm = \Lambda_0${\rm )} then ${\mathbb C}H^n_s (\lambda )$
consists of the two connected components $(\Lambda_+ \setminus
\Lambda_0 )/G_\lambda \approx S^{2n-1}_{2s} \times S^1$ and
$(\Lambda_- \setminus \Lambda_0 )/G_\lambda \approx
H^{2n-1}_{2s-1} \times S^1$. \label{t:2}
\end{theorem}
\noindent If ${\mathbb R}^N_\nu = ({\mathbb R}^N , h_{\nu , N})$,
with $h_{\nu , N} (x,y) = - \sum_{j=1}^\nu x_j y_j + \sum_{j=\nu
+1}^N x_j y_j$, then $S^N_\nu (r) = \{ x \in {\mathbb R}^{N+1} :
h_{\nu , N+1}(x,x) = r^2 \}$ ($r > 0$) is the pseudosphere in
${\mathbb R}^{N+1}_\nu$, while $H^N_\nu (r) = \{ x \in {\mathbb
R}^{N+1} : h_{\nu + 1, N+1} (x,x) = - r^2 \}$ ($r > 0$) is the
pseudohyperbolic space in ${\mathbb R}^{N+1}_{\nu +1}$. When $r =
1$ we write simply $S^N_\nu$ and $H^N_\nu$. Also $\odot$ denotes
the symmetric tensor product, e.g. $\alpha \odot \beta =
\frac{1}{2} (\alpha \otimes \beta + \beta \otimes \alpha )$ for
any $1$-forms $\alpha$ and $\beta$. A construction similar to that
in Theorem \ref{t:2} was performed in \cite{kn:DrFr} (for metrics
which are locally conformal to anti-K\"ahlerian metics, cf. Lemma
3, {\em op. cit.}, p. 119).
\par {\em Proof of Theorem} \ref{t:2}. If $F_\lambda^m (z) = z$
for some $z \in {\mathbb C}^n \setminus \Lambda_0$ then $m = 0$,
hence $G_\lambda$ acts freely on ${\mathbb C}^n \setminus
\Lambda_0$. Given $z_0 \in {\mathbb C}^n \setminus \Lambda_0$ let
$B_r (z_0 )$ be the open Euclidean ball of center $z_0$ and radius
$r$. Also, we set $\Omega_r (z_0 ) = B_r (z_0 ) \setminus
\Lambda_0$. As well known (cf. \cite{kn:KoNo}, Vol. II, p. 137)
$G_\lambda$ acts on ${\mathbb C}^n \setminus \{ 0 \}$ as a
properly discontinuous group of holomorphic transformations, hence
there is $r > 0$ such that $F^m_\lambda (B_r (z_0 )) \cap B_r (z_0
) = \emptyset$, and then $F^m_\lambda (\Omega_r (z_0 )) \cap
\Omega_r (z_0 ) = \emptyset$, for any $m \in {\mathbb Z} \setminus
\{ 0 \}$. As ${\mathbb C}^n \setminus \Lambda_0$ is an open subset
of ${\mathbb C}^n$ it follows (cf. e.g. \cite{kn:Boo}, p. 97) that
${\mathbb C}H^n_s (\lambda ) := ({\mathbb C}^n \setminus \Lambda_0
)/G_\lambda$ is a complex manifold.
\par
Assume from now on that $0 < \lambda < 1$. Let $\pi : {\mathbb
C}^n \setminus \Lambda_0 \to {\mathbb C}H^n_s (\lambda )$ be the
projection. To prove the last statement in Theorem \ref{t:2} we
consider the $C^\infty$ diffeomorphism
\[ F : {\mathbb C}H^n_s (\lambda ) \to \Sigma^{2n-1} \times S^1 ,
\]
\begin{equation}\label{e:diff}
F(\pi (z)) = \left( |z|_{s,n}^{-1} \, z \; , \; \exp \left(
\frac{2\pi i \log |z|_{s,n}}{\log \lambda} \right) \right) ,
\end{equation}
with the obvious inverse
\[ F^{-1}(\zeta , w) = \pi \left( \lambda^{\arg (w)/(2\pi )} \,
\zeta \right) , \;\;\; \zeta \in \Sigma^{2n-1} , \;\; w \in S^1 ,
\]
where $\arg : {\mathbb C} \to [0, 2 \pi )$. Finally, note that
$\Sigma^{2n-1} \cap \Lambda_+ = S^{2n-1}_{2s}$ and $\Sigma^{2n-1}
\cap \Lambda_- = H^{2n-1}_{2s-1}$. $\square$
\begin{proposition} The Lee form of $({\mathbb C}H^n_s (\lambda ),
g_{s,n})$ is locally given by
\begin{equation}
\omega = - d \log |z|^2_{s,n} . \label{e:11}
\end{equation}
In particular ${\mathbb C}H^n_s (\lambda )$ has a parallel Lee
form. Let $\Omega_{\pm} = (\Lambda_\pm \setminus \Lambda_0
)/G_\lambda$ be the connected components of ${\mathbb C}H^n_s
(\lambda )$ and $a(z) = {\rm sign}(b_{s,n}(z,z)) = \pm 1$ for $z
\in \Lambda_\pm \setminus \Lambda_0$. Then the Lee field $B$ of
$({\mathbb C}H^n_s (\lambda ) , g_{s,n})$ is given by
\begin{equation} B =  - 2 a(z) \left( z^j \frac{\partial}{\partial
z^j} + \overline{z}^j \frac{\partial}{\partial \overline{z}^j}
\right) . \label{e:12}
\end{equation}
Finally, if $\left. B_\pm = B\right|_{\Omega_\pm}$ then $B_+$ is
spacelike while $B_-$ is timelike. \label{p:1}
\end{proposition}
\noindent {\em Proof}.  An inspection of (\ref{e:10}) leads to
(\ref{e:11}) and hence to
\[ \omega = b_{s,n}(z,z)^{-1} \left\{ \sum_{j=1}^s (\overline{z}_j d z^j +
z_j d \overline{z}^j ) - \sum_{j= s+1}^n (\overline{z}_j d z^j +
z_j d \overline{z}^j ) \right\} , \] with the convention $z_j =
z^j$. Next
\[ g_{j\overline{k}} = \frac{1}{2} \, |z|_{s,n}^{-2} \,
\epsilon_j\,  \delta_{jk} \] (where $\epsilon_j = -1$ for $1 \leq
j \leq s$ and $\epsilon_j = 1$ for $s+1 \leq j \leq n$) yields
(\ref{e:12}). Since \[ Z_j \left( |z|^{-2}_{s,n} \right) = - a(z)
|z|_{s,n}^{-4} \, \epsilon_j \, \overline{z}_j \] the
identity
\[ 2 g_{AD} \Gamma_{BC}^A  = Z_B (g_{CD}) + Z_C (g_{BD}) - Z_D
(g_{BC}) \] leads to
\[ \Gamma^\ell_{jk} = - \frac{a(z)}{2 |z|^2_{s,n}}
\left( \epsilon_j \overline{z}_j \delta_k^\ell + \epsilon_k
\overline{z}_k \delta_j^\ell \right) , \;\;\;
\Gamma^{\overline{\ell}}_{jk} = 0, \]
\[ \Gamma_{j\overline{k}}^\ell =
\frac{a(z)}{2 |z|^2_{s,n}} \left( \epsilon_j \delta_{jk} z^\ell -
\epsilon_k z_k \delta_j^\ell \right) , \;\;
\Gamma^{\overline{\ell}}_{j\overline{k}} =  \frac{a(z)}{2
|z|^2_{s,n}} \left( \epsilon_j \delta_{jk} \overline{z}^\ell -
\epsilon_j \overline{z}_j \delta^\ell_k \right) , \] hence a
calculation shows that $\nabla_{Z_j} B = 0$. Finally $g_{s,n}(B,B)
= 4 a(z)$. $\square$

\vskip 0.1in Given two semi-Riemannian manifolds $M$ and $N$ and a
$C^\infty$ submersion $\Pi : M \to N$, we say $\Pi$ is a {\em
semi-Riemannian submersion} if 1) $\Pi^{-1} (y)$ is a
semi-Riemannian submanifold of $M$ for each $y \in N$, and 2) $d_x
\Pi : {\mathcal H}_x \to T_{\Pi (x)} (N)$ is a linear isometry of
semi-Euclidean spaces, where ${\mathcal H}_x = {\rm Ker}(d_x \Pi
)^\bot$, for any $x \in M$ (cf. \cite{kn:O'Ne}, p. 212). Also we
recall (cf. \cite{kn:BaRo}) the indefinite complex projective
space ${\mathbb C}P^{n-1}_s (k)$. Its underlying complex manifold
is the open subset of the complex projective space
\[ {\mathbb C}P^{n-1}_s (k)  = (\Lambda_+ \setminus \Lambda_0
)/{\mathbb C}^* \subset {\mathbb C}P^{n-1} \;\;\;\;\; ({\mathbb
C}^* = {\mathbb C} \setminus \{ 0 \} ). \] As to the
semi-Riemannian metric of ${\mathbb C}P^{n-1}_s (k)$, let
\[ \Pi : S^{2n-1}_{2s} ( 2/\sqrt{k}) \to
{\mathbb C}P^{n-1}_s (k), \;\;\; \Pi (z) = z \cdot {\mathbb C}^*
\;\;\; (k > 0)
\]
be the {\em indefinite Hopf fibration}. It is a principal
$S^1$-bundle and $S^1$ acts on $S^{2n-1}_{2s}(2/\sqrt{k})$ as a
group of isometries, hence (by slightly adapting the proof of
Proposition E.3 in \cite{kn:BGM}, p. 7, to the semi-Riemannian
context) there is a unique semi-Riemannian metric of index $2s$ on
${\mathbb C}P^{n-1}_s (k)$ such that $\Pi$ is a semi-Riemannian
submersion and ${\mathbb C}P^{n-1}_s (k)$ is an indefinite complex
space form of (constant) holomorphic sectional curvature $k$.
Again by \cite{kn:BaRo}, p. 57, ${\mathbb C}P^{n-1}_s (k)$ is
homotopy equivalent to ${\mathbb C}P^{n-1-s}$, hence ${\mathbb
C}P^{n-1}_s (k)$ is simply connected. We shall prove
\begin{theorem} Let $D = \{ 2\pi i a + (\log \lambda ) b : a , b
\in {\mathbb Z} \}$ {\rm (}$0 < \lambda < 1${\rm )} and consider
the torus $T^1_{\mathbb C} = {\mathbb C}/D$. Then $T^1_{\mathbb
C}$ acts freely on ${\mathbb C}H^n_s (\lambda )$ and
\[ p : \Omega_+ \to {\mathbb C}P^{n-1}_s (4), \;\;\; p(\pi (z)) =
z \cdot {\mathbb C}^* , \] is a principal $T^1_{\mathbb C}$-bundle
and a semi-Riemannian submersion of $\Omega_+$ {\rm (}carrying the
indefinite l.c.K. metric $g_{s,n}${\rm )} onto ${\mathbb
C}P^{n-1}_s (4)$. Moreover the complex Hopf manifold ${\mathbb
C}H^{n-s}(\lambda )$ {\rm (}respectively ${\mathbb C}H^s (\lambda
)${\rm )} is a strong deformation retract of $\Omega_+$ {\rm
(}respectively of $\Omega_-${\rm )} hence
\[ H_k (\Omega_+ ; {\mathbb Z}) = \begin{cases} {\mathbb Z}
\otimes {\mathbb Z}, & k = 2(n-s), \cr {\mathbb Z}, & k \neq
2(n-s), \cr
\end{cases}
\] \[ H_k (\Omega_- ; {\mathbb Z}) = \begin{cases} {\mathbb Z}
\otimes {\mathbb Z} , & k = 2s, \cr {\mathbb Z}, & k \neq 2s, \cr
\end{cases}
\] and $\Omega_{\pm}$ are not simply connected
\[ \pi_1 (\Omega_+ ) = \begin{cases} {\mathbb Z} \oplus {\mathbb
Z}, & s = n-1, \cr {\mathbb Z}, & s \neq n-1, \cr \end{cases}
\;\;\; \pi_1 (\Omega_- ) = \begin{cases} {\mathbb Z} \oplus
{\mathbb Z}, & s = 1, \cr {\mathbb Z}, & s \neq 1. \cr \end{cases}
\] Each fibre $p^{-1}(z \cdot {\mathbb C}^* )$, $z \in S^{2n-1}_{2s}$, is tangent
to the Lee field $B$ of $\Omega_+$ hence $p : \Omega_+ \to
{\mathbb C}P^{n-1}_s (4)$ is a harmonic map.
\label{t:3}
\end{theorem}
\noindent {\em Proof}. The action of $T^1_{\mathbb C}$ on
${\mathbb C}H^n_s (\lambda )$ is given by
\[ \pi (z) \cdot (\zeta + D) = \pi \left( e^\zeta \; z \right) ,
\;\;\; z \in {\mathbb C}^n \setminus \Lambda_0 , \;\; \zeta \in
{\mathbb C} . \] As $b_{s,n}(e^\zeta z , e^\zeta z) = e^{2 {\rm
Re}(\zeta )} b_{s,n}(z,z) \neq 0$ the action is well defined. To
see that the action is free, let us assume that $\pi (e^\zeta z_0
) = \pi (z_0 )$, for some $z_0 \in {\mathbb C}^n \setminus
\Lambda_0$. Then $e^\zeta _0 = \lambda^m z_0$, for some $m \in
{\mathbb Z}$, hence $\zeta = m \log \lambda + 2k \pi i$, for some
$k \in {\mathbb Z}$, i.e. $\zeta + D = 0$.
\par
To see that $S^1 \to \Omega_+ \stackrel{p}{\to} {\mathbb
C}^{n-1}_s (4)$ is a principal bundle let us assume that $p(\pi (z
)) = p( \pi (z^\prime ))$, with $z, z^\prime \in {\mathbb C}^n
\setminus \Lambda_0$. Then $z^\prime = \alpha z$, for some $\alpha
\in {\mathbb C}^*$. We wish to show that there is a unique $\zeta
+ D \in T^1_{\mathbb C}$ such that $\pi (z^\prime ) = \pi (z)
\cdot (\zeta + D)$. Indeed we may consider $\zeta = \log |\alpha |
+ i \arg (\alpha )$.
\par
Let $F_t : \Lambda_+ \setminus \Lambda_0 \to {\mathbb C}^n$, $0
\leq t \leq 1$, be given by
\[ F_t (z) = ((1-t) z^\prime \, , \, z^{\prime\prime}), \;\;\; z =
(z^\prime , z^{\prime\prime}) \in \Lambda_+ \setminus \Lambda_0 ,
\]
where $z^\prime = (z_1 , \cdots , z_s )$ and $z^{\prime\prime} =
(z_{s+1} , \cdots , z_n )$. Then
\[ b_{s,n}(F_t (z), F_t (z)) = - (1-t)^2 |z^\prime |^2 +
|z^{\prime\prime}|^2 \geq b_{s,n}(z,z) > 0 \] hence $F_t$ is
$(\Lambda_+ \setminus \Lambda_0 )$-valued. Therefore $F_t$ induces
a homotopy
\[ H^+_t : \Omega_+ \to \Omega_+ , \;\;\; H^+_t (\pi (z)) = \pi \left(
F_t (z) \right) , \;\;\; 0 \leq t \leq 1. \] Let us consider (cf.
e.g. \cite{kn:KoNo}, Vol. II, p. 137) the complex Hopf manifold
${\mathbb C}H^n (\lambda ) = ({\mathbb C}^n \setminus \{ 0 \}
)/G_\lambda$ and denote by $\pi_0 : {\mathbb C}^n \setminus \{ 0
\} \to {\mathbb C}H^n (\lambda )$ the projection. Let ${\mathbb
C}H^{n-s}(\lambda )$ be thought of as identified to
\[ \{ \pi_0 (z) \in
{\mathbb C}H^n (\lambda ) : z_1 = 0, \cdots , z_s = 0 \} . \] Note
that ${\mathbb C}H^{n-s}(\lambda ) \subset \Omega_+$. Also
\[ H^+_0 = 1_{\Omega_+} \, , \;\;\;
H^+_1 \left( \Omega_+ \right) \subset {\mathbb C}H^{n-s}(\lambda
),
\]
\[ \left. H^+_t \right|_{{\mathbb C}H^{n-s}(\lambda )} = i, \;\;\; 0 \leq t
\leq 1, \] (where $i : {\mathbb C}H^{n-s} (\lambda ) \to \Omega_+$
is the inclusion) hence ${\mathbb C}H^{n-s}(\lambda )$ is a strong
deformation retract of $\Omega_+$. Also $H^+_1 \circ i =
1_{{\mathbb C}H^{n-s}(\lambda )}$ and $H^+ : 1_{\Omega_+} \simeq i
\circ H^+_1$ (i.e. the maps $1_{{\mathbb C}H^{n-s}(\lambda )}$ and
$i \circ H^+_1$ are homotopic) so that $i, H^+_1$ are reciprocal
homotopy equivalences, i.e.
\begin{equation}
\Omega_+ \simeq {\mathbb C}H^{n-s}(\lambda ),
\label{e:13}
\end{equation}
(a homotopy equivalence). As well known, (\ref{e:13}) implies that
\[ i_* : H_k ({\mathbb C}H^{n-s}(\lambda ) ; {\mathbb Z}) \approx
H_k (\Omega_+ ; {\mathbb Z}), \] (a group isomorphism). Therefore,
to compute $H_k (\Omega_+ ; {\mathbb Z})$ it suffices to compute
the singular homology of the complex Hopf manifold. This is an
easy exercise in algebraic topology (based on the K\"unneth
formula). Indeed
\[ H_k ({\mathbb C}H^n (\lambda ) ; {\mathbb Z}) =
\sum_{p+q=k} H_p (S^{2n-1}; {\mathbb Z}) \otimes H_q (S^1 ;
{\mathbb Z}) =
\] \[ = H_{k-1}(S^{2n-1} ; {\mathbb Z}) \otimes H_1 (S^1 ; {\mathbb
Z}) = \begin{cases} {\mathbb Z} \otimes {\mathbb Z} , & k = 2n,
\cr {\mathbb Z}, & k \neq 2n, \cr \end{cases} \] yielding
(\ref{e:11}). As to the homotopy groups, again by (\ref{e:13})
\[ \pi_k (\Omega_+ ) \approx \pi_k ({\mathbb C}H^{n-s}) \]
(a group isomorphism) and if $n > 1$
\[ \pi_k ({\mathbb C}H^n (\lambda )) = \pi_k (S^{2n-1}) \oplus
\pi_k (S^1 ) = \begin{cases} {\mathbb Z}, & k \in \{ 1 , 2n-1 \} ,
\cr 0, & k \not\in \{ 1 , 2n-1 \} , \cr \end{cases} \] while if $n
= 1$
\[ \pi_k ({\mathbb C}H^1 (\lambda )) = \pi_k (S^1 ) \oplus \pi_k
(S^1 ) = \begin{cases} {\mathbb Z} \oplus {\mathbb Z}, & k = 1,
\cr 0, & k \neq 1. \cr \end{cases} \] Let us show now that $p :
\Omega_+ \to {\mathbb C}P^{n-1}_s (4)$ is a semi-Riemannian
submersion. Let $i : S^{2n-1}_{2s} \to \Lambda_+ \setminus
\Lambda_0$ be the inclusion and $V_{z_0} = {\rm Ker}(d_{z_0} \Pi
)$ the vertical space, $z_0 \in S^{2n-1}_{2s}$. As $\Pi :
S^{2n-1}_{2s} \to {\mathbb C}P^{n-1}_s (4)$ is a semi-Riemannian
submersion $V_{z_0}$ is nondegenerate, hence the perp space
$H_{z_0}$ of $V_{z_0}$ is also nondegenerate and
\[ T_{z_0} (S^{2n-1}_{2s}) = H_{z_0} \oplus_{\rm orth} V_{z_0} . \]
Let $V_{0, \pi (z_0 )} = {\rm Ker}(d_{\pi (z_0 )} p)$ and let
$N(S^{2n-1}_{2s})_{z_0}$ be the normal space of $i$ at $z_0$ (as
$S^{2n-1}_{2s}$ has index $2s$ it follows that
$N(S^{2n-1}_{2s})_{z_0}$ has index zero).
\begin{lemma} For any $z_0 \in S^{2n-1}_{2s}$
\[ V_{0, \pi (z_0 )} = (d_{z_0} \pi ) \{ N(S^{2n-1}_{2s})_{z_0}
\oplus (d_{z_0} i ) V_{z_0} \} . \] \label{l:2}
\end{lemma}
\noindent {\em Proof}. Let $N = z^j Z_j + \overline{z}^j
\overline{Z}_j$ be the unit normal on $S^{2n-1}_{2s}$ in
$\Lambda_+ \setminus \Lambda_0$, with the flat indefinite K\"ahler
metric
\[
g_0 = \sum_{j=1}^n \epsilon_j d z^j \odot d \overline{z}^j . \]
Let $U$ be the tangent vector field on $S^{2n-1}_{2s}$ defined by
$(d i) U = - J N$, where $J$ is the complex structure on ${\mathbb
C}^n$. Let
\[ \Pi_0 : \Lambda_+ \setminus \Lambda_0 \to {\mathbb
C}P^{n-1}_s (4) \]
be the canonical projection (so that $\Pi =
\left. \Pi_0 \right|_{S^{2n-1}_{2s}}$). Then $(d \Pi_0 ) U = 0$
hence $V_{z_0} = {\mathbb R} U_{z_0}$. In particular, by the
commutativity of the diagram
\[ \begin{array}{ccccc} {\mathbb C}H^n_s (\lambda ) & \supset &
\;\;\;\;\; \Omega_+ & \stackrel{p}{\longrightarrow} & {\mathbb
C}P^{n-1}_s (4)
\\
\empty & \empty & \pi \uparrow & \empty & \uparrow \Pi \\
{\mathbb C}^n_s & \supset & \;\;\; \Lambda_+ \setminus \Lambda_0 &
\stackrel{i}{\longleftarrow} & S^{2n-1}_{2s} \end{array} \] it
follows that
\begin{equation} (d_{z_0} \pi )(d_{z_0} i) V_{z_0} \subseteq V_{0,
\pi (z_0 )} . \label{e:14}
\end{equation}
On the other hand $\Pi_0$ is a holomorphic map hence
\[ (d_{\pi (z_0 )} p)(d_{z_0} \pi ) N_{z_0} = (d_{z_0} \Pi_0 )
J_{z_0} (d_{z_0} i) U_{z_0} = \]
\[ = J^\prime_{\Pi_0 (z_0 )} (d_{z_0} \Pi_0 )(d_{z_0} i) U_{z_0} =
J^\prime_{\Pi (z_0 )} (d_{z_0} \Pi )U_{z_0} = 0, \] where
$J^\prime$ denotes the complex structure on ${\mathbb C}P^{n-1}$.
We obtain
\begin{equation}
(d_{z_0} \pi ) N(S^{2n-1}_{2s})_{z_0} \subseteq V_{0, \pi (z_0 )}
. \label{e:15}
\end{equation}
At this point Lemma \ref{l:2} follows from
(\ref{e:14})-(\ref{e:15}) and an inspection of dimensions.
$\square$
\par
As $|z_0 |_{s,t} = 1$ the indefinite scalar product $(\; , \;
)_{z_0}$ induced on $T_{z_0} (\Lambda_+ \setminus \Lambda_0 )$ by
$g_{0, z_0}$ coincides with that induced by $|z|_{s,t}^{-2} \sum_j
\epsilon_j d z^j \odot d \overline{z}^j$ at $z_0$.
\begin{lemma} $V_{0, \pi (z)}$ is nondegenerate, for any $z \in \Lambda_+ \setminus
\Lambda_0$. \label{l:3}
\end{lemma}
\noindent {\em Proof}. Let $A = - J B$, where $J$ is the complex
structure on $\Omega_+$. Then $\{ A_{\pi (z)} , B_{\pi (z)} \}$
span $V_{0, \pi (z)}$ and (by Proposition \ref{p:1}) both $A$ and
$B$ are spacelike. $\square$
\par
By Lemma \ref{l:3} the perp space $H_{0, \pi (z_0 )}$ of $V_{0,
\pi (z_0 )}$ is also nondegenerate and
\[ T_{\pi (z_0 )} (\Omega_+ ) = H_{0, \pi (z_0 )} \oplus_{\rm
orth} V_{0, \pi (z_0 )} . \] Let $v \in H_{z_0}$. Then $(d_{z_0}
i)v$ is perpendicular on $N(S^{2n-1}_{2s})_{z_0} \oplus (d_{z_0} i
) V_{z_0}$. On the other hand $d_{z_0} \pi : T_{z_0}(\Lambda_+
\setminus \Lambda_0 ) \to T_{\pi (z_0 )} (\Omega_+ )$ is a linear
isometry hence (by Lemma \ref{l:2}) $(d_{z_0} \pi )(d_{z_0} i )v$
is perpendicular on $V_{0, \pi (z_0 )}$, and then it lies on
$H_{0, \pi (z_0 )}$. Again by inspecting dimensions we obtain
\begin{equation}
H_{0, \pi (z_0 )} = (d_{z_0} \pi )(d_{z_0} i) H_{z_0} .
\label{e:16}
\end{equation}
Next, by (\ref{e:16}) and by $d_{z_0} \Pi : H_{z_0} \approx T_{\Pi
(z_0 )} ({\mathbb C}P^{n-1}_s (4))$ (a linear isometry), it
follows that
\begin{equation}\label{e:17} d_{\pi (z )} p : H_{0, \pi (z )} \approx T_{\Pi (z )}
({\mathbb C}P^{n-1}_s (4)), \end{equation} a linear isometry for
any $z \in S^{2n-1}_s$. We wish to show that (\ref{e:17}) actually
holds for any $z \in \Lambda_+ \setminus \Lambda_0$.
\begin{lemma} The torus $T^1_{\mathbb C}$ acts on $\Omega_+$ as a group of
isometries of the semi-Riemannian manifold $(\Omega_+ , g_{s,n})$.
\label{l:4}
\end{lemma}
\noindent {\em Proof}. If $g = w + D \in T^1_{\mathbb C}$ ($w \in
{\mathbb C}$) then the right translation $R_g : \Omega_+ \to
\Omega_+$ is a holomorphic map locally given by $z \mapsto z e^w$
hence
\[ |z e^w |^{-2}_{s,n} \; \left( (d R_g ) Z_j \; , \; (d R_g ) \overline{Z}_k
\right)_{z e^w} = \] \[ = e^{2 {\rm Re}(w)} \; |z|^{-2}_{s,n} \;
b_{s,n} (e^w e_j , e^w e_k ) = |z|^{-2}_{s,n} \; \epsilon_j \;
\delta_{jk}
\] where $\{ e_j : 1 \leq j \leq n \}$ is the canonical basis in ${\mathbb C}^n$.
$\square$
\par
As $p \circ R_g =$ const., each $R_g$ preserves the vertical
spaces $V_0$. Then (by Lemma \ref{l:4}) $R_g$ preserves the
horizontal spaces $H_0$, as well. Therefore, to complete the proof
we must show that for any $z \in \Lambda_+ \setminus \Lambda_0$
there is $g \in T^1_{\mathbb C}$ and $z_0 \in S^{2n-1}_{2s}$ such
that $\pi (z) = \pi (z_0 ) g$. Indeed we may consider $g = \log
|z|_{s,t} + D$ and $z_0 = |z|_{s,t}^{-1} z$.
\par
To prove the last statement in Theorem \ref{t:2} we establish
\begin{lemma} Let $z_0 \in \Lambda_+ \setminus \Lambda_0$
and let $j : T^1_{\mathbb C} \to \Omega_+$ be the immersion given
by $j (\zeta + D) = \pi (z_0 ) \cdot (\zeta + D)$. Then
\[ \left. (d j) \frac{\partial}{\partial u} \right|_{\zeta + D} = -
\frac{1}{4}\; e^\zeta \; B_{\pi (z_0 )} \, , \] where $\zeta = u +
i v$, hence $j(T^1_{\mathbb C})$ is tangent to the Lee field of
$\Omega_+$.
\end{lemma} This follows easily from (\ref{e:12}). Then (by
Proposition \ref{p:3}) $j : T^1_{\mathbb C} \to \Omega_+$ is a
minimal isometric immersion. Therefore $p$ is a semi-Riemannian
submersion with minimal fibres, hence a harmonic map (in the sense
of \cite{kn:Fug}). Compare to Theorem 3 in \cite{kn:Dra}, p. 375.
$\square$

\section{An indefinite l.c.K. metric with nonparallel Lee form}
Let ${\mathbb C}_+ = \{ w \in {\mathbb C} : {\rm Im}(w) > 0 \}$ be
the upper half space and consider the {\em Tricerri metric} (cf.
\cite{kn:DrOr}, p. 24) \[ g_{0,1} = {\rm Im}(w)^{-2} \; d w \odot
d \overline{w} + {\rm Im}(w) \; d z \odot d \overline{z} . \]
$g_{0,1}$ is (by a result in \cite{kn:Tri}) a (positive definite)
g.c.K. metric on ${\mathbb C}_+ \times {\mathbb C}$. We build a
family of indefinite l.c.K. metrics of index $0 \leq s < n$
containing the Tricerri metric as a limiting case (for $s = 0$ and
$n=1$).
\begin{proposition} Let $0 \leq s < n$ and $n \geq 1$. Let
$g_{s,n}$ be the indefinite Hermitian metric on ${\mathbb C}_+
\times {\mathbb C}^n_s$ given by
\begin{equation} g_{s,n} = {\rm Im}(w)^{-2} \; d w \odot d \overline{w} +
{\rm Im}(w) \sum_{j=1}^n \epsilon_j d z^j \odot d \overline{z}^j .
\label{e:19}
\end{equation}
Then $g_{s,n}$ is an indefinite globally conformal K\"ahler metric
with a nonparallel Lee form and its Lee field is spacelike.
Moreover, let $a \in {\rm SL}(3 , {\mathbb Z})$ be a unimodular
matrix with ${\rm Spec}(a) = \{ \alpha , \beta , \overline{\beta}
\}$, where $\alpha > 1$ and $\beta \in {\mathbb C} \setminus
{\mathbb R}$. Let $G_{\alpha , \beta}$ be the group of holomorphic
transformations of ${\mathbb C}_+ \times {\mathbb C}^n_s$
generated by $F_0 (w , z) = (\alpha w , \beta z)$, $w \in {\mathbb
C}_+$, $z \in {\mathbb C}^n$. Then $g_{s,n}$ is $G_{\alpha ,
\beta}$-invariant. \label{p:2}
\end{proposition}
\noindent {\em Proof}. Note that we may write (\ref{e:19}) as
$g_{s,n} = {\rm Im}(w) g_0$ where
\[ g_0 = {\rm Im}(w)^{-3} \; d w \odot d \overline{w} +
\sum_{j=1}^n \epsilon_j d z^j \odot d \overline{z}^j \] is an
indefinite Hermitian metric whose K\"ahler $2$-form $\Omega_0$ is
\[ - i \{ {\rm Im}(w)^{-3} \; d w \wedge d \overline{w} +
\sum_{j=1}^n \epsilon_j \; d z^j \wedge d \overline{z}^j \} .
\] Hence $d \Omega_0 = 0$, i.e. $g_0$ is an indefinite K\"ahler
metric. Thus (\ref{e:19}) is an indefinite l.c.K. metric whose Lee
form
\[ \omega = d f = \frac{1}{w - \overline{w}} \, (d w - d
\overline{w}) \;\;\;\; (f = \log {\rm Im}(w)) \] is exact. Raising
indices we obtain the Lee field
\[ B = i \, {\rm Im}(w) \, \left( \frac{\partial}{\partial w} -
\frac{\partial}{\partial \overline{w}} \right) \] so that
$g_{s,n}(B,B) = 1$, i.e. $B$ is spacelike. Next, the only
surviving coefficients of the Levi-Civita connection $\nabla$ of
$({\mathbb C}_+ \times {\mathbb C}^n_s , g_{s,n})$ are
\[ \Gamma_{k0}^j = - \Gamma_{k\overline{0}}^j = - \frac{i}{4} \,
{\rm Im}(w)^{-1} \, \delta_k^j \, , \] hence $\nabla_{Z_j} B =
\frac{1}{2} \, Z_j \neq 0$. Finally, the $G_{\alpha ,
\beta}$-invariance of $g_{s,n}$ follows from $F_0^* d w = \alpha
\; d w$, $F^*_0 d z^j = \beta \; d z^j$, $1 \leq j \leq n$, and
from $\alpha \beta \overline{\beta} = \det (a) = 1$. $\square$

\section{The second canonical foliation}
Let $(M , J, g)$ be an indefinite l.c.K. manifold, of index $\nu =
2 s$. Let $A = - J B$ be the {\em anti-Lee field}. We also set
$\theta (X) = g(X , A)$ (the {\em anti-Lee form}), so that $\theta
= \omega \circ J$. Also, let $\Omega (X,Y) = g(X , J Y)$ be the
K\"ahler $2$-form. Since $D J = 0$ it follows that
\begin{equation}
(\nabla_X J) Y = \frac{1}{2} \{ \theta (Y) X - \omega (Y) J X -
g(X,Y) A - \Omega (X,Y) B \} \label{e:20}
\end{equation}
for any $X,Y \in T(M)$. As an immediate application of
(\ref{e:20}) we have
\begin{proposition} Let $(M , J , g)$ be an indefinite l.c.K. manifold and
$i : N \hookrightarrow M$ a complex submanifold of $M$ such that
$i^* g$ is a semi-Rie\-man\-nian metric. Let $h$ be the second
fundamental form of $i$. Then \begin{equation} h(J X , J Y) = -
h(X,Y) - g(X,Y) B^\bot \, , \label{e:18}
\end{equation}
for any $X,Y \in T(N)$, where $B^\bot$ is the normal component of
the Lee field of $M$. Then the mean curvature vector of $i$ is
given by $H = - \frac{1}{2} B^\bot$. In particular $i$ is minimal
if and only if $N$ is tangent to the Lee field. \label{p:3}
\end{proposition} \noindent
This extends a result of I. Vaisman (cf. \cite{kn:Vai2}) to the
case of semi-Riemannian complex submanifolds of an indefinite
l.c.K. manifold. {\em Proof of Proposition} \ref{p:3}. As $N$ is a
complex manifold $T(N)$ admits a local orthonormal frame of the
form $\{ E_\alpha , J E_\alpha : 1 \leq \alpha \leq m \}$, i.e.
$g(E_\alpha , E_\beta ) = \epsilon_\alpha \delta_{\alpha\beta}$,
and then (by (\ref{e:18})) the mean curvature vector $H$ of $i$ is
given by
\[ H = \frac{1}{2m} \; \sum_\alpha \epsilon_\alpha \{ h(E_\alpha , E_\alpha ) + h(J
E_\alpha , J E_\alpha ) \} = - \frac{1}{2} \; B^\bot . \] It
remains that we prove (\ref{e:18}). Let ${\rm tan}$ and ${\rm
nor}$ be the projections associated with the decomposition $T(M) =
T(N) \oplus T(N)^\bot$ and let us set
\[ t \xi = {\rm tan}(\xi ), \;\;\; f \xi = {\rm nor}(\xi ), \;\;\;
\xi \in T(N)^\bot . \] Then $J \xi = t \xi + f \xi$ and by
applying $J$ once more we get $f^2 = - I$. Let $A^\bot = {\rm
nor}(A)$ and $B^\bot = {\rm nor}(B)$. Then (by (\ref{e:20}) and
the Gauss formula)
\[ h(X , J Y) = f \; h(X,Y) - \frac{1}{2} \{ g(X, Y) A^\bot +
\Omega (X,Y) B^\bot \} , \] for any $X,Y \in T(N)$. Finally, using
$A^\bot = - f \; B^\bot$ and $f^2 = - I$ we obtain (\ref{e:18}).
$\square$ \vskip 0.1in \noindent As another application of
(\ref{e:20}) we shall prove
\begin{theorem}
Let $M$ be a complex $n$-dimensional {\rm (}$n > 1${\rm )}
indefinite l.c.K. manifold with a parallel Lee form $(\nabla
\omega = 0)$ and $c = g(B,B) \in {\mathbb R}$. Let $M(c) = M
\setminus {\rm Sing}(B)$, an open subset of $M$. Then
\[ P : x \in M(c) \mapsto {\mathbb R} A_x \oplus {\mathbb R}B_x
\subset T_x (M) \] is an integrable distribution, hence $P$
determines a foliation ${\mathcal F}_c$ of $M(c)$ by real surfaces
such that {\rm i)} either $c \neq 0$ and then each leaf $L \in
M/{\mathcal F}_c$ is Riemannian {\rm (}with the metric ${\rm
sign}(c) \; i^* g$, $i : L \hookrightarrow M${\rm )} and a totally
geodesic surface in $(M , g)$, or {\rm ii)} $c = 0$ and then each
leaf $L \in M(c)/{\mathcal F}_c$ is either an isotropic surface
{\rm (}when $n \geq 3${\rm )} or a totally lightlike surface {\rm
(}when $n=2${\rm )}. Assume that $n \geq 3$. Then the second
fundamental form of a leaf $L \in M(0)/{\mathcal F}_0$ with
respect to any transversal vector bundle ${\rm tr}(T(L)) \to L$
vanishes. \label{t:4last}
\end{theorem}
\noindent {\em Proof}. Note that $c \neq 0$ yields $M(c) = M$.
Moreover (by the very definition of the anti-Lee field) $\{ A_x ,
B_x \}$ are linearly dependent if and only if $x \in {\rm
Sing}(B)$. Hence the sum ${\mathbb R} A_x + {\mathbb R}B_x$ is
direct, for any $x \in M(c)$. To see that $P$ is involutive it
suffices to check that $[A, B] \in P$. Let $X \in P^\bot$. Then
\[ g([A,B] , X) = g(\nabla_A B - \nabla_B A , X) = \;\;\;\;\;\;\;
({\rm as} \; \nabla B = 0, \; \nabla g = 0) \]
\[ = - B(g(A,X)) + g(A , \nabla_B X ) = \;\;\;\;\;\;\;\;\;\;\;\;
({\rm as} \; \theta (X) = 0) \]
\[ = - g(J B , \nabla_B X ) = g(B , J \nabla_B X) = \;\;\;\;\;\;\;\;\;\;\;\;
\;\;\;\;\;\;\;\;\;\;\;\; ({\rm by} \; (\ref{e:20})) \]
\[ = g(B , \nabla_B J X ) = B(g(B , J X)) - g(\nabla_B B , J X) =
B(\theta (X)) = 0, \] hence $[A,B] \in (P^\bot )^\bot = P$. By the
classical Frobenius theorem there is a foliation ${\mathcal F}_c$
of $M(c)$ such that $P = T({\mathcal F}_c )$. Assume that $c \neq
0$. Then either $P$ is spacelike (when $c > 0$) or timelike (when
$c < 0$). Let $L \in M/{\mathcal F}_c$ and let $h$ be the second
fundamental form of $L \hookrightarrow M$. Then $\nabla B = 0$
yields $h(A,B) = h(B,B) = 0$. Finally (by (\ref{e:20}))
\[ \nabla_A A = - \nabla_A J B = \] \[ = - J \nabla_A B - \frac{1}{2} \{
\theta (B) A - \omega (B) J A - g(A,B) A - \Omega (A,B) B \} = 0
\] so that $h(A,A) = 0$. We may conclude that $h = 0$. $\square$
\par Assume now that $c = 0$, so that both the Lee and anti-Lee
fields are lightlike. Let us set
\[ {\rm Rad} \; P_x = P_x \cap P^\bot_x \, , \;\;\; x \in M(0). \]
We have ${\rm dim}_{\mathbb R} P_x = 2$ and ${\rm dim}_{\mathbb R}
P^\bot_x = 2(n-1)$, hence ${\rm Rad} \; P = P$. Therefore each
leaf $L \in M(0)/{\mathcal F}_0$ is a $2$-lightlike submanifold
(surface) in $M(0)$ and in particular (according to the
terminology in \cite{kn:BeDu}, p. 149-150) an isotropic
submanifold (when $n
> 2$) or a totally lightlike submanifold (when $n = 2$) of $M(0)$.
We shall need the following adaptation of Lemma 1.4 and Theorem
1.6 in \cite{kn:BeDu}, p. 149-150 (to the case of the lightlike
foliation ${\mathcal F}_0$, rather than a single istropic
submanifold)
\begin{lemma} Assume that $n \geq 3$.
Let $S(P^\bot ) \to M(0)$ be a vector subbundle of
$P^\bot \to M(0)$ such that
\begin{equation}
P^\bot = P \oplus_{\rm orth} S(P^\bot ). \label{e:21}
\end{equation}
Then for any $x \in M(0)$ there exist an open neighborhood $U
\subseteq M(0)$ and a system of linearly independent tangent
vector fields $\{ N_1 , N_2 \}$ on $U$ such that
\begin{equation}
\theta (N_1 ) = \omega (N_2 ) = 1, \;\;\; \theta (N_2 ) = \omega
(N_1 ) = 0, \label{e:22}
\end{equation}
\begin{equation}
g(N_i , N_j ) = 0, \;\;\; g(N_i , W) = 0, \label{e:23}
\end{equation}
for any $W \in S(P^\bot )$. Moreover, if $\left. {\rm
ltr}(P)\right|_U$ is given by
\[ {\rm ltr}(P)_x = {\mathbb R}N_{1, x} \oplus {\mathbb R} N_{2,
x} \, , \;\;\; x \in U, \] then the vector bundles $\left. {\rm
ltr}(P)\right|_U$ glue up to a vector bundle ${\rm ltr}(P) \to
M(0)$ such that
\begin{equation}
S(P^\bot )^\bot = P \oplus {\rm ltr}(P). \label{e:24}
\end{equation}
\label{l:6}
\end{lemma}
\noindent According to the terminology in \cite{kn:BeDu}, $\left.
{\rm ltr}(P)\right|_L$ is the {\em lightlike transversal vector
bundle} of $(L , \left. S(P^\bot )\right|_L )$, for each leaf $L
\in M(0)/{\mathcal F}_0$. By Proposition 2.1 in \cite{kn:BeDu}, p.
5, $S(P^\bot )$ is nondegenerate. Let $S(P^\bot )^\bot$ be the
orthogonal complement of $S(P^\bot )$. Of course, this is also
nondegenerate and
\begin{equation}
\label{e:25} T(M(0)) = S(P^\bot ) \oplus_{\rm orth} S(P^\bot
)^\bot .
\end{equation}
Note that $S(P^\bot )^\bot$ has rank $4$ and (by (\ref{e:21})) $P
\subset S(P^\bot )^\bot$. To prove Lemma \ref{l:6} let $E \to
M(0)$ be a subbundle of $S(P^\bot )^\bot \to M(0)$ such that
\begin{equation}
S(P^\bot )^\bot = P \oplus E. \label{e:26}
\end{equation}
For any $x \in M(0)$ there is an open neighborhood $U \subseteq
M(0)$ of $x$ and a local frame $\{ V_1 , V_2 \} \subset
\Gamma^\infty (U, E)$. Let $D \in C^\infty (U)$ be given by
\[ D = \theta (V_1 ) \omega (V_2 ) - \omega (V_1 ) \theta (V_2 ). \]
We claim that $D(x) \neq 0$, for any $x \in U$. Indeed, if $D(x_0
) = 0$ for some $x_0 \in U$ then
\[ \theta (V_2 )_{x_0} = \lambda \theta (V_1 )_{x_0} \, , \;\;\;
\omega (V_2 )_{x_0} = \lambda \omega (V_1 )_{x_0} \, , \] for some
$\lambda \in {\mathbb R}$, hence $v_\lambda := V_{2, x_0} -
\lambda V_{1, x_0}$ is orthogonal to both the Lee and anti-Lee
vectors, i.e. $v_\lambda \in P^\bot_{x_0}$. Yet $v_\lambda \in
E_{x_0} \subset S(P^\bot )^\bot_{x_0}$, i.e. $v_\lambda$ is
orthogonal to $S(P^\bot )_{x_0}$. Then (by (\ref{e:21}))
$v_\lambda \in P_{x_0} \cap E_{x_0} = (0)$, i.e. $\{ V_{1, x_0} ,
V_{2 , x_0} \}$ are linearly dependent, a contradiction. Let us
set
\[ N_1 = \lambda_{11} A + \lambda_{12} B + \frac{1}{D} \{ \omega
(V_2 ) V_1 - \omega (V_1 ) V_2 \} , \]
\[ N_2 = \lambda_{21} A + \lambda_{22} B - \frac{1}{D} \{ \theta
(V_2 ) V_1 - \theta (V_1 ) V_2 \} , \] where $\lambda_{ij} \in
C^\infty (U)$ are given by
\[ \lambda_{11} = - \frac{1}{D^2} \{ \omega (V_2 )^2 g(V_1 , V_1 )
- 2 \omega (V_1 ) \omega (V_2 ) g(V_1 , V_2 ) + \omega (V_1 )^2
g(V_2 , V_2 ) \} , \]
\[ \lambda_{22} = - \frac{1}{D^2} \{ \theta (V_2 )^2 g(V_1 , V_1 )
- 2 \theta (V_1 ) \theta (V_2 ) g(V_1 , V_2 ) + \theta (V_1 )^2
g(V_2 , V_2 ) \} , \]
\[ \lambda_{12} = \lambda_{21} = \frac{1}{2 D^2} \{ \omega (V_2 )
\theta (V_2 ) g(V_1 , V_1 ) + \] \[ + [\theta (V_1 ) \omega (V_2 )
+ \omega (V_1 ) \theta (V_2 ) ] g(V_1 , V_2 ) - \omega (V_1 )
\theta (V_1 ) g(V_2 , V_2 ) \} . \] A calculation shows that $N_i$
are linearly independent at each $x \in U$ and satisfy
(\ref{e:22})-(\ref{e:23}). Therefore $\left. {\rm
ltr}(P)\right|_U$ is well defined. Let $U^\prime \subseteq M(0)$
be another open neighborhood of $x$ and $\{ V_1^\prime ,
V_2^\prime \}$ a local frame of $E$ on $U^\prime$, so that
$V^\prime_i = f_i^j V_j$, for some $f^j_i \in C^\infty (U \cap
U^\prime )$. A calculation shows that $\lambda^\prime_{ij} =
\lambda_{ij}$ and then $N^\prime_i = N_i$ on $U \cap U^\prime$,
hence $\left. {\rm ltr}(P)\right|_U$ and $\left. {\rm
ltr}(P)\right|_{U^\prime}$ glue up over $U \cap U^\prime$.
Finally, one may check that $P_{x_0} \cap {\rm ltr}(P)_{x_0} \neq
(0)$ at some $x_0 \in U$ yields $D(x_0 ) = 0$, a contradiction.
Hence the sum $P + {\rm ltr}(P)$ is direct and (\ref{e:24}) must
hold. $\square$ \vskip 0.1in With the choices in Lemma \ref{l:6}
we set
\[ {\rm tr}(P) = {\rm ltr}(P) \oplus_{\rm orth} S(P^\bot ), \]
so that (cf. \cite{kn:BeDu1}) $\left. {\rm tr}(P)\right|_L$ is the
{\em transversal bundle} of $L \in M(0)/{\mathcal F}_0$. Then
(\ref{e:25}) yields $T(M(0)) = P \oplus {\rm tr}(P)$ and we may
decompose $\nabla_X Y = \nabla^P_X Y + h^P (X,Y)$, for any $X,Y
\in P$, such that $\nabla^P$ is a connection in $P \to M(0)$ and
$h^P$ is a $C^\infty (M(0))$-bilinear symmetric ${\rm
tr}(P)$-valued form on $P$ (compare to (2.1) in \cite{kn:BeDu1},
p. 154). Once again we may use $\nabla B = 0$, $\nabla_A A = 0$ to
conclude that $h^P = 0$. $\square$

\section{A CR extension result}
Let $M$ be a complex $n$-dimensional indefinite l.c.K. manifold of
index $2s$, $0 < s < n$, with a parallel Lee form. Let $\mathcal
F$ be the first canonical foliation (given by $\omega = 0$). Each
leaf $L$ of $\mathcal F$ is a real hypersurface in $M$, hence a CR
manifold with the CR structure \[ T_{1,0}(L) = T^{1,0}(M) \cap
[T(L) \otimes {\mathbb C}] \] induced by the complex structure of
$M$. There is a natural first order differential operator
\[ \overline{\partial}_L : C^1 (L) \to \Gamma^\infty (T_{0,1}(L)^*)
\]
given by $(\overline{\partial}_L f) \overline{Z} =
\overline{Z}(f)$, for any $C^1$ function $f : L \to {\mathbb C}$
and any $Z \in T_{1,0}(L)$. Here $T_{0,1}(L) =
\overline{T_{1,0}(L)}$. The solutions to $\overline{\partial}_L f
= 0$ are the CR functions on $L$ (and $\overline{\partial}_L f =
0$ are the tangential Cauchy-Riemann equations on $L$, cf. e.g.
\cite{kn:Bog}, p. 124). Let $CR^k (L)$ be the space of all CR
functions on $L$, of class $C^k$. It is a natural question whether
a CR function on a leaf $L$ of $\mathcal F$ extends to a
holomorphic function on $M$ (at least locally). We answer this
question for the canonical foliation of $\Omega_+$ (a similar
result holds for $\Omega_-$) where an explicit description of the
leaves is available. Precisely
\begin{theorem} Let $n \geq 2$ and $0 < s < n$ such that $s \neq (n-1)/2$ when
$n$ is odd. Let $w \in S^1$ and $N_w^+ = F^{-1} (S^{2n-1}_{2s}
\times \{ w \} )$, where $F : \Omega_+ \to S^{2n-1}_{2s} \times
S^1$ is the diffeomorphism {\rm (\ref{e:diff})}. Let $\mathcal F$
be the foliation of $\Omega_+$ given by the Pfaff equation $d \log
|z|^2_{s,n} = 0$. Then the leaf space is
\begin{equation}
\label{e:28}
\Omega_+ /{\mathcal F} = \{ N^+_w : w \in S^1 \}
\end{equation}
and for any point $x \in N_w^+$ there is an open neighborhood $U$
of $x$ in $M$ such that for any $f \in CR^1 (N_w^+ )$ there is a
holomorphic function $F \in {\mathcal O}(U)$ such that $\left. F
\right|_{U \cap N_w^+} = f$. \label{t:5}
\end{theorem}
\noindent {\em Proof}. Note that
\[ N^+_w = \{ \pi (\lambda^{{\rm arg}(w)/(2\pi )} \zeta ) : \zeta
\in S^{2n-1}_{2s} \} \;\;\; (w \in S^1 ). \] Let $z \in \Lambda_+
\setminus \Lambda_0$ and let us consider $\zeta = |z|^{-1}_{s,t} z
\in S^{2n-1}_{2s}$ and $w = \exp (2\pi i \log |z|_{s,n} / \log
\lambda ) \in S^1$. Then $\arg (w) = 2\pi \log |z|_{s,n} /\log
\lambda + 2 m \pi$ for some $m \in {\mathbb Z}$, so that
\[ \pi (z) = \pi (\lambda^m |z|_{s,t} \zeta ) = \pi
(\lambda^{\arg (w)/(2 \pi )} \zeta ) \in N_w^+ \, , \] that is
through each point $\pi (z) \in \Omega_+$ passes at least one
hypersurface of the form $N_w^+$. Next, let us assume that $\pi
(z) \in N^+_w \cap N_{w^\prime}^+$. Then
\[ e^{\arg (w^\prime )/(2\pi )} \zeta^\prime = \lambda^m e^{\arg
(w)/(2 \pi )} \zeta \] for some $\zeta , \zeta^\prime \in
S^{2n-1}_{2s}$ and some $m \in {\mathbb Z}$. Then $b_{s,n}(\zeta ,
\zeta ) = b_{s,n}(\zeta^\prime , \zeta^\prime ) = 1$ imply
\[ \arg (w^\prime ) = \arg (w) + 2 m \pi \log \lambda  \]
hence $N_w^+ = N^+_{w^\prime}$, that is through each $\pi (z ) \in
\Omega_+$ passes a unique hypersurface of the form $N_w^+$. To
emphasize, $N^+_w = N^+_{w^\prime}$ if and only if $w^\prime =
e^{2m\pi i \log \lambda} w$, for some $m \in {\mathbb Z}$.
Therefore, to prove (\ref{e:28}) it suffices to check that the Lee
field $B$ of $\Omega_+$ is orthogonal to each $N_w^+$. We set
\[ D(0, r) = \{ z \in \Lambda_+ \setminus \Lambda_0 : |z|_{s,n} <
r \} \;\;\;\; (r > 0) \] and consider the annulus $A_k = D(0 ,
\lambda^k ) \setminus \overline{D}(0, \lambda^{k+1})$, $k \in
{\mathbb Z}$. If $U_k = \pi (A_k )$ then $\phi_k = (\pi : A_k \to
U_k )^{-1}$ are local charts on $\Omega_+$. Note that the
holomorphic transformation $F_\lambda$ maps the pseudosphere
$S^{2n-1}_{2s}(\lambda^k )$ onto $S^{2n-1}_{2s}(\lambda^{k+1})$,
for any $k \in {\mathbb Z}$. In other words, when building
$\Omega_+$ one identifies the points where a generic complex line
through the origin intersects the pseudospheres
$S^{2n-1}_{2s}(\lambda^k )$. In particular $\pi
(S^{2n-1}_{2s}(\lambda^k )) = \pi (S^{2n-1}_{2s})$ and $U_k =
U_0$, for any $k \in {\mathbb Z}$.
\begin{lemma} Let $w \in S^1$ and $a = \arg (w)/(2\pi \log \lambda
)$. If $a \in {\mathbb R} \setminus {\mathbb Z}$ then $N_w^+
\subset U_0$, while if $a \in {\mathbb Z}$ then $N_w^+ = \pi
(S^{2n-1}_{2s})$. In particular, for any $w \in S^1 \setminus \{
e^{2m\pi i \log \lambda} : m \in {\mathbb Z} \}$
\begin{equation}
\label{e:phi0}
\phi_0 (N_w^+ ) = S^{2n-1}_{2s}\left( \lambda^{- [a]} e^{\arg
(w)/(2\pi )} \right) .
\end{equation}
\label{l:7}
\end{lemma}
\noindent Here $[a]$ is the integer part of $a \in {\mathbb R}$.
Note that (\ref{e:phi0}) doesn't apply to the leaf $L_0 = \pi
(S^{2n-1}_{2s}) \in \Omega_+ /{\mathcal F}$ (corresponding to $a
\in {\mathbb Z}$). However, in this case one may consider $\lambda
< \epsilon < 1$ and the annulus $A = D(0, \lambda^{-1} \epsilon )
\setminus \overline{D}(0, \epsilon )$ and then $L_0$ is contained
in $U = \pi (A)$ and $(\pi : A \to U)^{-1}$ is a local chart on
$\Omega_+$. To prove Lemma \ref{l:7} let $x$ be a point of
$N_w^+$, $x = \pi (e^{\arg (w)/(2\pi )} \zeta )$, and let us set
$z = \lambda^{k - [a]} e^{\arg (w)/(2 \pi )} \zeta$. Then $a - 1 <
[a] \leq a$ yields $\lambda^{k+1} < |z|_{s,n} \leq \lambda^k$ and
the second inequality becomes an equality if and only if $a \in
{\mathbb Z}$, that is if $w = e^{2m\pi i \log \lambda}$, for some
$m \in {\mathbb Z}$. Finally, if $a \in {\mathbb R} \setminus
{\mathbb Z}$ then $N_w^+ \subset U_0$ and
\[ \phi_0 (N^+_w ) = \{ \phi_0 (\pi (^{\arg (w)/(2 \pi )} \zeta )) :
\zeta \in S^{2n-1}_{2s} \} = \] \[ = \{ \lambda^{-[a]} e^{\arg
(w)/(2 \pi )} \zeta : \zeta \in S^{2n-1}_{2s} \}  =
S^{2n-1}_{2s}\left( \lambda^{-[a]} e^{\arg (w)/(2\pi )} \right)
\] and $B = - 2(z^j Z_j + \overline{z}^j \overline{Z}_j )$ is
orthogonal to any $S^{2n-1}_{2s}(r)$. The Cayley transform
\[ {\mathcal C}(z) = \left( \frac{z^\prime}{r + z_n} \, , \,
\frac{i(r - z_n )}{r + z_n} \right) , \;\;\; z = (z^\prime , z_n )
\in {\mathbb C}^n \setminus \{ z_n + r = 0 \}, \] is a CR
isomorphism of $S^{2n-1}_{2s}(r)$ onto $\partial {\mathcal
S}_{s,n} \setminus \{ \zeta_n + i = 0 \}$, where
\[ {\mathcal S}_{s,n} = \{ \zeta \in {\mathbb C}^n : {\rm
Im}(\zeta_n ) > \sum_{\alpha =1}^{n-1} \epsilon_\alpha
|\zeta_\alpha |^2 \} . \] The CR structure $T_{1,0}(\partial
{\mathcal S}_{s,n})$ is the span of $\{ \partial /\partial
\zeta^\alpha + 2 i \epsilon_\alpha \overline{\zeta}_\alpha
\partial /\partial \zeta^n : 1 \leq \alpha \leq n - 1 \}$ hence the Levi
form has signature $(s , n-s-1)$. Yet $s \geq 1$ hence (by H.
Lewy's CR extension theorem, cf. e.g. Theorem 1 in \cite{kn:Bog},
p. 198) for any $x \in \partial {\mathcal S}_{s,n}$ there is an
open neighborhood $U \subseteq {\mathbb C}^n$ of $x$ such that for
any $f \in CR^1 (\partial {\mathcal S}_{s,n})$ there is a unique
$F \in {\mathcal O}(U \cap {\mathcal S}_{s,n}) \cap C^0 (U \cap
\overline{\mathcal S}_{s,n})$ such that $\left. F\right|_{U \cap
\partial {\mathcal S}_{s,n}} = f$. Then the last statement in Theorem \ref{t:5}
holds for any $x \in N_w^+ \setminus \{ \pi (e^{\arg (w)/(2\pi )}
\zeta ) : \zeta \in \Lambda^{n-1}_0 \times \{ - 1 \} \}$, where
$\Lambda^{n-1}_0 = \Lambda^{n-1} \cup \{ 0 \}$ and $\Lambda^{n-1}$
is the null cone in ${\mathbb C}^{n-1}_s$ (so that $\phi_0 (x)$
satisfies $z_n + r \neq 0$ ($r = \lambda^{-[a]} e^{\arg (w)/(2\pi
)}$)). For arbitrary $x \in N_w^+$ the argument requires that $z_j
+ r \neq 0$, for some $1 \leq j \leq n$ (the remaining case is
ruled out by our assumption that $n \neq 2s+1$). $\square$

\section{Levi foliations}
Let $(M , J , g)$ be a complex $n$-dimensional indefinite l.c.K.
manifold and $B$, $A$ its Lee and anti-Lee fields, respectively.
Let us set $Z := B + i A \in T^{1,0}(M)$. Clearly $\omega (Z ) =
c$. Let us assume that $\nabla \omega = 0$ and ${\rm Sing}(B) =
\emptyset$ and set
\[ T_{1,0}({\mathcal F}) = T^{1,0}(M) \cap [T({\mathcal F})
\otimes {\mathbb C}] \] so that the portion of $T_{1,0}({\mathcal
F})$ over a leaf $L \in M/{\mathcal F}$ is the CR structure of
$L$. Also the portion of $H({\mathcal F}) := {\rm Re}\{
T_{1,0}({\mathcal F}) \oplus \overline{T_{1,0}({\mathcal F})} \}$
over $L$ is the Levi distribution $H(L)$ of $L$. The distribution
$H({\mathcal F})$ carries the complex structure
\[ J : H({\mathcal F}) \to H({\mathcal F}), \;\;\; J(V +
\overline{V}) = i(V - \overline{V}), \;\;\; V \in
T_{1,0}({\mathcal F}). \] See also \cite{kn:DrNi}. Let us set
\[ {\mathcal L}(V , \overline{W}) = i \; \pi \; [V ,
\overline{W}], \;\;\; V , W \in T_{1,0}({\mathcal F}), \] where
$\pi : T({\mathcal F}) \to T({\mathcal F}) /H({\mathcal F})$ is
the natural projection, so that ${\mathcal L}$ is the Levi form of
each leaf of $\mathcal F$. The {\em null space} of $\mathcal L$ is
\[ {\rm Null}({\mathcal L}) = \{ V \in T_{1,0}({\mathcal F}) : {\mathcal L}(V , \overline{V}) = 0 \} . \]
We may state the following corollary of Theorems \ref{t:1} and
\ref{t:4last}
\begin{proposition} If the Lee vector $B$ is lightlike then the Levi
form of each leaf of $\mathcal F$ is degenerate $(Z \in {\rm
Null}({\mathcal L}))$ and ${\mathcal F}_0$ is a subfoliation of
$\mathcal F$. Moreover if $n = 2$ then each leaf of $\mathcal F$
is Levi-flat and the Levi foliation of each leaf $L \in
M/{\mathcal F}$ extends to a unique holomorphic foliation of $M$.
\label{p:4}
\end{proposition}
We recall that a CR manifold $L$ is {\em Levi-flat} if its Levi
form vanishes identically (${\mathcal L} = 0$). If this is the
case $L$ is foliated by complex manifolds (whose complex dimension
equals the CR dimension of $L$). The resulting foliation (whose
tangent bundle if the Levi distribution $H(L)$ of $L$) is the {\em
Levi foliation} of $L$. If $L$ is embedded in some come complex
manifold $M$ a problem raised by C. Rea, \cite{kn:Rea}, is whether
the Levi foliation of (a Levi-flat CR manifold) $L$ may extend to
a holomorphic foliation of $M$. Proposition \ref{p:4} exhibits a
family of Levi foliations of class $C^\infty$ which extend
holomorphically (while Rea's extension theorem (cf. {\em op.
cit.}) requires real analytic data).
\par
{\em Proof of Proposition} \ref{p:4}. Let us assume that $c = 0$.
Then $Z \in T({\mathcal F}) \otimes {\mathbb C}$ hence $Z \in
T_{1,0}({\mathcal F})$. Moreover $Z + \overline{Z} \in H({\mathcal
F})$ yields $B \in H({\mathcal F})$ and by applying $J$ we may
conclude that $A \in H({\mathcal F})$ as well. Hence (with the
notations of Theorem \ref{t:4last}) $P \subseteq H({\mathcal F})$.
Note that $[Z , \overline{Z}] = 2 i [A,B]$ and then (by the
integrability of $P$) ${\mathcal L}(Z , \overline{Z}) = 0$. When
$n = 2$ each leaf $L$ of $\mathcal F$ is a $3$-dimensional CR
manifold hence $H({\mathcal F}) = P$ and $L$ is Levi-flat. Finally
the foliation induced by ${\mathcal F}_0$ on $L$ is precisely the
Levi foliation of $L$. In other words, the Levi foliation of each
leaf of $\mathcal F$ extends to a holomorphic foliation of $M$
which is precisely the second canonical foliation of $M$.

{\small {\bf Acknowledgements}. The first author was partially
supported by INdAM (Italy) within the interdisciplinary project
{\em Nonlinear subelliptic equations of variational origin in
contact geometry}. The second author was supported by Natural
Sciences and Engineering Council of Canada. This paper was
completed while the first author was a guest of the Department of
Mathematics of the University of Windsor (Ontario, Canada).}

\end{document}